\documentclass[12pt]{amsart}
\usepackage{amsfonts,latexsym}

\newtheorem{theorem}{Theorem}
\newtheorem{lemma}[theorem]{Lemma}
\newtheorem{thm}[theorem]{Theorem}

\newtheorem{proposition}[theorem]{Proposition} 
 
\newtheorem{corollary}[theorem]{Corollary}

\theoremstyle{definition}
\newtheorem{definition}[theorem]{Definition}

\theoremstyle{remark}
\newtheorem{remark}[theorem]{Remark}

\newcommand{\brfr}{$\hspace{0 pt}$} 

\DeclareMathOperator{\cf}{cf}

\title[Combinatorial principles, compactness of spaces IV]
{Combinatorial and model-theoretical principles related to 
regularity of ultrafilters and compactness of topological spaces. IV.}

\author[]{Paolo Lipparini} 
\address{Dipartimento di Matematica\\
Viale della Ricerca Scientifica\\
II Uniservitaccia di Roma (Tor Vergata)\\
I-00133 ROME 
ITALY
}
\urladdr{http://www.mat.uniroma2.it/\textasciitilde lipparin}
\thanks{The author has received support from MPI and GNSAGA.
We wish to express our gratitude to X. Caicedo for stimulating discussions and correspondence} 


\begin{document} 

\begin{abstract} 
We extend to singular cardinals the model-theoretical relation 
$\lambda \stackrel{\kappa}{\Rightarrow} \mu$
introduced in \cite{bumi}. We extend some results obtained in 
Part II, finding equivalent conditions involving uniformity of ultrafilters and the existence of
 certain
 infinite matrices. Our present definition suggests a new
compactness property for abstract logics. 
\end{abstract}

\maketitle




See Parts I, II, III \cite{parti} or 
\cite{CN,CK,KM,BF,easter,bumi,arxiv}
 for unexplained notation.

For $ \lambda $, $\nu$ infinite cardinals,
$S_\nu( \lambda )$ denotes
the set of subsets of $ \lambda $ of cardinality $< \nu$.
A set $X \subseteq S_\nu( \lambda )$ is \emph{cofinal} in 
$ S_\nu( \lambda )$      
if and only if for every $x \in S_\nu( \lambda )$
there exists $y \in X$ such that $x \subseteq y$.
The \emph{cofinality} $\cf S_\nu( \lambda )$ of $ S_\nu( \lambda )$
is the minimal cardinality of a cofinal subset of $ S_\nu( \lambda )$.
Notice that if $ \lambda $ is regular, then
$\cf S_ \lambda  ( \lambda )= \lambda $.

\begin{remark}\label{unif} 
If $X \subseteq S_\lambda ( \lambda )$ is cofinal in 
$ S_\lambda ( \lambda )$ then an ultrafilter $D$ over $ \lambda $ is uniform 
if and only if $ x \not \in D$ for every $x \in X$, if and only if 
$ \lambda \setminus x \in D$ for every $x \in X$.
\end{remark}

We are now going to extend to singular cardinals the
definition of $\lambda \stackrel{\kappa}{\Rightarrow} \mu$
as given in \cite{bumi} and recalled in Part II (see also \cite{easter},
 \cite[Section 0]{arxiv}).

\begin{definition}\label{limpm} 
Suppose that 
$ \lambda \geq \mu$ are infinite cardinals, and 
$ \kappa \geq \sup \{ \cf S_ \lambda  ( \lambda ), \cf S_\mu(\mu)\} $.

Fix some set $V \subseteq S_ \lambda  ( \lambda )$
 cofinal in  $ S_ \lambda  ( \lambda )$
of cardinality $\leq \kappa $, 
and some set $W \subseteq S_ \mu ( \mu )$ cofinal in  $ S_ \mu ( \mu )$
of cardinality $\leq \kappa $. 
(Of course, it is always possible to choose $V$ of cardinality
$\cf S_ \lambda  ( \lambda )$ and $W$ of cardinality $ \cf S_\mu(\mu)$.)

Consider the model 
$ \mathfrak A'= \langle  \lambda, R, R_v \rangle _{v \in V \cup W}   $,
where

(i) for $ \gamma < \lambda $,
$ \mathfrak A' \models R( \gamma )$
if and only if $ \gamma < \mu$,
and 

(ii) for $ \gamma < \lambda $ and $ v \in V \cup W $,
$ \mathfrak A' \models R_v( \gamma )$
if and only if $ \gamma \in v$. 

Now for the definition of the relation $\lambda \stackrel{\kappa}{\Rightarrow} \mu$.

The relation $\lambda \stackrel{\kappa}{\Rightarrow} \mu$
 means that the model $ \mathfrak A'$  has an
expansion ${\mathfrak A}$ in a language with at most $ \kappa $ new symbols such that whenever
$\mathfrak B $
is an elementary extension of
$\mathfrak A$ and $ \mathfrak B$ has an element $x$ such that 
$ {\mathfrak B} \models \neg R_v( x) $ for every $ v \in V$, 
then   
$ \mathfrak B$ has an element $y$ such that 
$ {\mathfrak B} \models R(y)$
and 
$ {\mathfrak B} \models \neg R_w(y) $ for every $ w \in W$.
\end{definition} 

\begin{remark}\label{independent}
Notice that the above definition is independent
from the choice of the sets $V$ and $W$.
Indeed, suppose that 
$\lambda \stackrel{\kappa}{\Rightarrow} \mu$
holds relative to some given particular choice of 
$V$ 
cofinal in  $ S_ \lambda  ( \lambda )$
of cardinality $\leq \kappa $, 
and  $W$ cofinal in  $ S_ \mu ( \mu )$
of cardinality $\leq \kappa $. Notice that such a choice is always
possible, since 
$ \kappa \geq \sup \{ \cf S_ \lambda  ( \lambda ), \cf S_\mu(\mu)\} $.

Consider another choice $V^*$, $W^*$ of cofinal sets of
cardinality $\leq \kappa $, and
consider the model 
$ \mathfrak C'= \langle  \lambda, R, R_v \rangle _{v \in V^* \cup W^*}$
such that, as above, 
$ \mathfrak C' \models R_v( \gamma )$
if and only if $ \gamma \in v$. 

Expand 
$ \mathfrak C'$ to  a model
$ \mathfrak C$ 
by adding (at most $ \kappa $-many new) relations $R_v$ ($v \in V \cup W$),
plus all constants, functions and relations 
added to $ \mathfrak A'$ in order to 
get some $ \mathfrak A$
witnessing 
$\lambda \stackrel{\kappa}{\Rightarrow} \mu$
relative to the choice of $V$, $W$ (notice that $ \mathfrak A'$ and
$ \mathfrak C'$ have the same base set.)   

Whenever $\mathfrak D $
is an elementary extension of
$\mathfrak C$, 
then an appropriate reduct $\mathfrak B $ of $\mathfrak D $
is an elementary extension of
$\mathfrak A$.
It is now easy to show that if   
$\mathfrak A$
witnesses 
$\lambda \stackrel{\kappa}{\Rightarrow} \mu$
when the definition 
is given relative to $V$, $W$,
then 
$\mathfrak C$
witnesses 
$\lambda \stackrel{\kappa}{\Rightarrow} \mu$
when the definition 
is given relative to $V^*$, $W^*$ 
(since $V$ and $W$ are cofinal).

Thus, the definition of
$\lambda \stackrel{\kappa}{\Rightarrow} \mu$
is independent from the choice of the pair of the cofinal sets,
provided, of course, that they have cardinality $\leq \kappa $.
\end{remark}

\begin{remark}\label{equivbumi}
We are now going to show that, when 
$\kappa \geq \lambda \geq \mu$, and 
both $ \lambda $ and $\mu$
are regular cardinals,
then the relation
$\lambda \stackrel{\kappa}{\Rightarrow} \mu$
as defined here is equivalent to
$\lambda \stackrel{\kappa}{\Rightarrow} \mu$
as defined in \cite[Definition 1.2]{bumi}
(assuming $ \kappa \geq \lambda $ is no essential loss of generality
since, in the sense of \cite[Definition 1.2]{bumi},
$\lambda \stackrel{\kappa}{\Rightarrow} \mu$
is equivalent to 
$\lambda \stackrel{\kappa'}{\Rightarrow} \mu$,
when $\omega \leq \kappa, \kappa' \leq \lambda $).

First observe that if $\lambda $ is regular, then
$\cf S_ \lambda  ( \lambda )= \lambda $, 
hence, if both $ \lambda $ and $\mu$
are regular, and $ \lambda \geq \mu$, then the condition  
$ \kappa \geq \sup \{ \cf S_ \lambda  ( \lambda ), \brfr \cf S_\mu(\mu)\} $
in Definition \ref{limpm} becomes simply
$ \kappa \geq \lambda $. 

Moreover, if $\lambda $ is regular, 
then the set $ V=\{ [0, \delta ]| \delta < \lambda \} $
is cofinal in 
$ S_ \lambda  ( \lambda )$;
similarly, if 
$\mu$ is regular, then 
$ W=\{ [0, \alpha ] | \alpha < \mu \} $
is  cofinal in $S_\mu(\mu)$.
Letting $\mathfrak A'$ be as in Definition \ref{limpm}
for the above choice of $V$ and $W$, we have
that if $v=[0, \delta ]$ then $\neg R_v(\gamma )$
holds in 
$\mathfrak A'$
if and only if 
$\delta < \gamma $.
Similarly, if $w \in W$, say,
$w=[0, \varepsilon ]$, then 
$R( \alpha  ) \wedge \neg R_w( \alpha )$
holds in 
$\mathfrak A'$
if and only if 
$\varepsilon  < \alpha  < \mu$.

Now notice that the definition of 
$\lambda \stackrel{\kappa}{\Rightarrow} \mu$
as given here involves 
an element $x$ such that 
$ {\mathfrak B} \models \neg R_v( x) $ for every $ v \in V$, 
and an element  
$y$ such that 
$ {\mathfrak B} \models R(y)$
and 
$ {\mathfrak B} \models \neg R_w(y) $ for every $ w \in W$, while
the definition of 
$\lambda \stackrel{\kappa}{\Rightarrow} \mu$
as given in \cite{bumi} (or in Part II) involves 
an element $x$ such that 
$ {\mathfrak B} \models \gamma < x $ for every $ \gamma < \lambda $, 
and an element $y$ such that 
$ {\mathfrak B} \models \alpha  < y < \mu $ for every $ \alpha  < \mu $.

Hence, the above comments show that, 
for the above particular choice of the sets
$V$ and $W$, and modulo appropriate expansions,
the elements $x$ and $y$ are required to satisfy 
sentences which are equivalent. Hence
the two definitions of 
$\lambda \stackrel{\kappa}{\Rightarrow} \mu$
are equivalent, since, as noticed in Remark \ref{independent}, our present definition 
does not depend on the choice of the sets $V$ and $W$.
\end{remark}

\begin{remark}\label{regsing}
On the contrary, for either $ \lambda  $ or  $\mu$ singular, the notion
$\lambda \stackrel{\kappa}{\Rightarrow} \mu$
introduced in Definition \ref{limpm}
is \textbf{not} equivalent to the notion 
\mbox{$(\lambda, \lambda ) \stackrel{\kappa}{\Rightarrow} (\mu, \mu)$}
introduced in \cite{easter} (see also \cite[Section 0]{arxiv}).

Indeed, $\lambda \stackrel{\kappa}{\Rightarrow} \mu$,
as introduced here,
involves analogues of $ \lambda $-\hspace{0 pt}decomposability of ultrafilters, 
while 
$(\lambda, \lambda ) \stackrel{\kappa}{\Rightarrow} (\mu, \mu)$
involves analogues of $( \lambda , \lambda )$-regularity.

Actually, when $ \kappa \geq \max \{2^ \lambda,  2^\mu\}   $, then 
$\lambda \stackrel{\kappa}{\Rightarrow} \mu$
is equivalent to 
``\emph{every $ \lambda $-decomposable ultrafilter is $ \mu$-decomposable}'' (e. g. by Theorem \ref{lmkprod4}), while
$(\lambda, \lambda ) \stackrel{\kappa}{\Rightarrow} (\mu, \mu)$
is equivalent to 
``\emph{every $( \lambda, \lambda ) $-regular ultrafilter is 
$ (\mu, \mu)$-regular}'' (\cite[Proposition 1]{easter}). Cf. also \cite[Theorem 2.5]{bumi}.

These notions are distinct for singular cardinals, since, for example,
every $( \lambda  ^+, \lambda  ^+)$-regular ultrafilter is 
 $( \lambda  , \lambda  )$-regular, while if $ \kappa $ 
is strongly compact, there exists a $ \kappa ^{ \omega +1} $-decomposable
ultrafilter which is not  $ \kappa ^{ \omega} $-decomposable.

However, it is not difficult to show that the two notions  
$\lambda \stackrel{\kappa}{\Rightarrow} \mu$
and
$(\lambda, \lambda ) \stackrel{\kappa}{\Rightarrow} (\mu, \mu)$
coincide
when both $ \lambda $ and $\mu$ 
are  regular cardinals. 
When 
$ \kappa \geq 2^ \lambda $
this is obvious, since, for $ \lambda $ regular,
an ultrafilter is $ \lambda $-decomposable if and only if 
it is  $( \lambda , \lambda )$-regular.
\end{remark}

\begin{remark}\label{6x}
Let $ \lambda $ be an infinite cardinal, and let $V \subseteq S_\lambda(\lambda )$ be
cofinal in  $ S_ \lambda  ( \lambda )$. Let
$ \mathfrak A'= \langle  \lambda,  R_v \rangle _{v \in V }   $,
where the $ R_v $'s are as in Definition 
\ref{limpm}.
It is easy to show that, 
an ultrafilter $D$ is $\lambda $-decomposable if and only if 
in the model $ \mathfrak B = \prod_D \mathfrak A'$ there is an element
$x$ such that $ {\mathfrak B} \models \neg R_v( x) $ for every $ v \in V$.
 \end{remark}

\begin{thm}\label{lmkprod4}
Suppose that  $ \lambda \geq \mu$  are infinite cardinals, 
and $ \kappa \geq \brfr \sup \{ \cf S_ \lambda  ( \lambda ), \cf S_\mu(\mu)\} $. 
Then the following conditions are equivalent.

\smallskip

(a) $\lambda \stackrel{\kappa}{\Rightarrow} \mu$ holds.

\smallskip

(b) There are $ \kappa $ functions $ (f_ \beta ) _{ \beta \in \kappa } $
from $ \lambda $ to $\mu$ such that whenever $D$ is an ultrafilter
uniform over $ \lambda $ then there exists some $ \beta \in \kappa $
such that $f_ \beta (D)$ is uniform over $ \mu$.

\smallskip

(b$'$) There are $ \kappa $ functions $ (f_ \beta ) _{ \beta \in \kappa } $
from $ \lambda $ to $\mu$ such that 
for every function $g: \kappa \to S_\mu(\mu)$ there exists some finite
set $F \subseteq \kappa $ such that 
$ \left| \bigcap _{\beta \in F} f_\beta  ^{-1}(g(\beta )) \right| < \lambda $.

\smallskip

(c) There is a family $ (C_{ \alpha , \beta }) _{ \alpha \in \mu , \beta \in \kappa}  $ 
of subsets of $ \lambda $ such that:

(i) For every $ \beta \in \kappa$, 
$ (C_{ \alpha , \beta }) _{ \alpha \in \mu }  $
is a partition of $ \lambda $.

(ii) For every function $g : \kappa  \to S_\mu (\mu) $ there exists a finite subset
$F \subseteq \kappa  $ such that 
$|\bigcap _{\beta \in F} \bigcup _{ \alpha \in g( \beta )} C_{ \alpha  , \beta }| < \lambda $.
\end{thm} 

\begin{proof}
Let us fix
$V \subseteq S_ \lambda  ( \lambda ) $ 
cofinal in  $ S_ \lambda  ( \lambda )$
of cardinality $\leq \kappa $, 
and  $W \subseteq S_ \mu ( \mu )$ cofinal in  $ S_ \mu ( \mu )$
of cardinality $\leq \kappa $.
Such sets $V$ and $W$ exist by the assumption 
$ \kappa \geq \sup \{ \cf S_ \lambda  ( \lambda ), \cf S_\mu(\mu)\} $.

Notice that, as we proved in Remark \ref{independent} , the definition of 
$\lambda \stackrel{\kappa}{\Rightarrow} \mu$
is independent from the choice of the sets $V$ and $W$.

(a) $\Rightarrow$ (b).
Let ${\mathfrak A}$
be an expansion of 
$  \langle  \lambda, R, R_v \rangle _{v \in V \cup W}   $
witnessing
$\lambda \stackrel{\kappa}{\Rightarrow} \mu$.

Without loss of generality (since $\kappa \geq \lambda$) we can assume that 
each element of $A$ is represented by some constant symbol, and that
${\mathfrak A}$ has Skolem functions (see \cite[Section 3.3]{CK}).
Indeed, since $ \kappa \geq |V \cup W|$, adding Skolem functions
to ${\mathfrak A}$ involves adding at most $ \kappa $ new symbols.

Consider the set of all functions $f: \lambda \to \mu$ 
which are definable in ${\mathfrak A}$. Enumerate them as 
$ (f_ \beta ) _{ \beta \in \kappa } $. We are going to show that 
these functions witness (b).

Indeed, let $D$ be an ultrafilter
uniform over $ \lambda $. 
Consider the $D$-class $Id_D$ of the identity
function  on $ \lambda $. 
Since $D$ is
uniform over $ \lambda $, letting
  $ {\mathfrak C} = \prod_D \mathfrak A$  we have that
$ \mathfrak C \models \neg R_v (Id_D)$ for every $ v \in V$
(since $|v| \leq \lambda $, for $v \in V$).

Let $ {\mathfrak B} $ be the Skolem hull
of $Id_D$ in ${\mathfrak C}$. 
By 
\L o\v s Theorem,
$ {\mathfrak C} \equiv {\mathfrak A} $.
Since 
$\mathfrak A $ has Skolem functions, 
hence $\mathfrak C $ has Skolem functions, then
$ {\mathfrak B} \equiv {\mathfrak C}$, by
 \cite[Proposition 3.3.2]{CK}. By transitivity,
$ {\mathfrak B} \equiv {\mathfrak A} $.
Since $ \mathfrak A$ has a name for each element of $A$,
then   $ {\mathfrak B} $ is an elementary
extension of $  {\mathfrak A} $.

Since ${\mathfrak A}$
witnesses
$\lambda \stackrel{\kappa}{\Rightarrow} \mu$,
and since 
$ \mathfrak B \models \neg R_v (Id_D)$ for every $ v \in V$
then   
$ \mathfrak B$ has an element $y_D$ such that 
$ {\mathfrak B} \models R(y_D) $
and 
$ {\mathfrak B} \models \neg R_w(y_D) $
for every $w \in W$.

Since $ {\mathfrak B} $ is the Skolem hull
of $Id_D$ in ${\mathfrak C}$, we have $y_D =f(Id_D)$,
that is, $y_D =f_D$, for
some function $f: \lambda \to A= \lambda $ definable in ${\mathfrak A}$.
Since $f$ is definable, then also the following function $f'$
is definable:
\[
f'( \gamma )=
\begin{cases}
f( \gamma) & \textrm{ if } f(\gamma )< \mu \\
0 & \textrm{ if } f(\gamma ) \geq \mu \\
\end{cases}
\]

Since $ {\mathfrak B} \models  R(y_D)$, then
$ \{ \gamma < \lambda | \mathfrak A \models R(y( \gamma ))\} \in D$,
that is,
\mbox{$ \{ \gamma < \lambda | y( \gamma )< \mu \} \in D$}.
Since
$y_D =f_D$,
$ \{ \gamma < \lambda | y( \gamma )=f( \gamma ) \} \in D$.
Hence,
$ \{ \gamma < \lambda | y( \gamma ) = f'( \gamma )\} \in D$,
since it contains the intersection of two sets in $D$.
Thus, 
$y_D =f'_D$.

Since $f': \lambda \to \mu$ and
$f'$ is definable  in ${\mathfrak A}$, then 
$f= f_ \beta $ for some $ \beta \in \kappa $, thus
$y_D =(f_ \beta)_D $.

We need to show that $D'= f_ \beta (D)$ is uniform over $\mu$.
Indeed, suppose that $ z \in S_ \mu(\mu)$. We have to show that
$z \not \in f_ \beta (D)$.

Since $W$ is supposed to be cofinal in $ S_ \mu(\mu)$,
there is $w \in W$ such that $z \subseteq w$.
Since $ {\mathfrak B} \models  \neg R_w(y_D)$, then
$ \{ \gamma < \lambda | \mathfrak A \models \neg R_w( y( \gamma )) \} \in D $;
that is,
$ \{ \gamma < \lambda | \mathfrak A \models \neg R_w( f_ \beta ( \gamma )) \} \in D $, that is, 
$ \{ \gamma < \lambda | f_ \beta ( \gamma ) \not \in w \} \in D $,
that is,
$ \{ \alpha  < \mu | \alpha  \not \in w \} \in f_ \beta (D) $,
that is, $w \not \in f_ \beta (D)$.
Since $z \subseteq w$, then 
$z \not \in D'=f_ \beta (D)$.

We have showed that no element of 
$S_\mu (\mu)$ belongs to $D'$
and this means
that $D'$ is uniform over $ \mu$.

(b) $\Rightarrow$ (a). Suppose we have functions
$ (f_ \beta ) _{ \beta \in \kappa } $ as given by (b).

Expand 
$  \langle  \lambda, R, R_v \rangle _{v \in V \cup W}   $
to a model 
$ \mathfrak A$ 
by adding, for each $ \beta \in \kappa $, a new function symbol representing
$f_ \beta $ (by abuse of notation, in what follows we shall
write $f_ \beta $ both for the function itself and
for the symbol that represents it).

Suppose that $\mathfrak B $ is an elementary extension of  
$ \mathfrak A$ and 
$ \mathfrak B$ has an element $x$ such that 
$ {\mathfrak B} \models \neg R_v( x) $ for every $ v \in V$.

For every formula $ \varphi (z)$ with just one variable $z$ 
in the language of  
$ \mathfrak A$ let 
$E_ \varphi = \{ \gamma < \lambda |  \mathfrak A \models \varphi ( \gamma )\}  $.
Let $F= \{ E_ \varphi | \mathfrak B \models \varphi (x)\} $.
Since the intersection of any two members of $F$ is still in 
$F$, and $\emptyset \not\in F$
(since $ \mathfrak A \equiv \mathfrak B$ and since if
$\mathfrak B \models \phi(x)$ then $\mathfrak B \models \exists z \phi(z)$),
 then $F$ can be extended to 
an ultrafilter $D$ on $ \lambda $. 

For every $ v \in V$,
consider the formula $ \varphi (z) \equiv \neg R_v( z)$. We get
$E_ \varphi = \{ \gamma < \lambda |  
\mathfrak A \models \neg R_v(  \gamma ) \}=
\lambda \setminus v $.
On the other side, 
since 
$ {\mathfrak B} \models \neg R_v( x) $, then
by the definition of $F$ we have
$ E_ \varphi = \lambda \setminus v \in F \subseteq D$.
Thus, by Remark \ref{unif}, and since $V$
is cofinal in $S_ \lambda (\lambda )$, $D$ is uniform over $ \lambda $.  

By (b), 
$ f_ \beta (D) $ is uniform over  $ \mu$, for
some $ \beta \in \kappa$. Thus,
$\mu \setminus w \in f_ \beta (D)$, 
for every $w \in W$.
That is,
$ \{ \gamma < \lambda | \mathfrak A \models \neg R_w (f_ \beta ( \gamma ))\} \in D $
for every $w \in W$.

For every $w \in W$,
consider the formula 
$ \psi (z) \equiv \neg R_w (f_ \beta ( z))$.
By the previous sentence, 
$E_ \psi \in D$.
Notice that 
$E _{\neg \psi}  $ is the complement of 
$E_ \psi $ in $ \lambda $.
Since $D$ is proper,
and $ E_ \psi \in D$, then
$E _{\neg \psi} \not\in D $.
Since $D$ extends $F$,
and 
either $E_ \psi \in F$
or $E _{\neg \psi}  \in F$,
 we necessarily have
 $E_ \psi \in F$, that is, 
$\mathfrak B \models \psi (x)$, that is, 
$\mathfrak B \models \neg R_w (f_ \beta ( x))$.

Since $ w \in W$ has been chosen
arbitrarily, we have that
$\mathfrak B \models \neg R_w (f_ \beta ( x))$
for every $ w \in W$.
Moreover, since $ f_ \beta : \lambda \to \mu$,
then $ \mathfrak A \models \forall z R(f_\beta (z))$,
hence, since $ \mathfrak B \equiv \mathfrak A$, then 
$\mathfrak B \models R(f_ \beta (x))$.

Thus, we have proved that 
$ \mathfrak B$ has an element $y= f_ \beta (x)$ such that 
$\mathfrak B \models R(y)$ and
$ {\mathfrak B} \models \neg R_w (y) $ for every $w \in W $.

(b) $\Leftrightarrow$ (b$'$) follows from Lemma \ref{lem} below.

(b$'$) $ \Rightarrow $ (c).
Suppose that we have functions $ (f_ \beta ) _{ \beta \in \kappa } $ as
given by (b$'$).
For $\alpha \in \mu $
and $\beta \in \kappa$, define 
$ C_{ \alpha , \beta } = f_ \beta ^{-1} ( \{ \alpha \} )  $. 

The family $ (C_{ \alpha , \beta }) _{ \alpha \in \mu , \beta \in \kappa}$
trivially satisfies Condition (i). Moreover, 
Condition (ii) is clearly equivalent to the condition
imposed on the $f_ \beta $'s in  (b$'$).

(c) $ \Rightarrow $ (b$'$).
Suppose we are given the family $ (C_{ \alpha , \beta }) _{ \alpha \in \mu , \beta \in \kappa}$
from (c). For $ \beta \in \kappa $ and $ \gamma < \lambda $,  define
$f_ \beta ( \gamma )$ to be 
the unique  (by (i)) $ \alpha \in \mu$  
such that $ \gamma \in C_{ \alpha , \beta }$.
Then 
$ f_\beta  ^{-1}(g(\beta ))= 
\bigcup _{ \alpha \in g( \beta )} C_{ \alpha  , \beta }$, thus 
Condition (ii) in (c)
implies that  (b$'$) holds.
\end{proof} 

\begin{corollary}\label{coruf}
If $ \lambda \geq \mu$ and $ \kappa \geq 2^ \lambda $,
then
$\lambda \stackrel{\kappa}{\Rightarrow} \mu$,
if and only if every ultrafilter uniform
over $ \lambda $ is $\mu$-decomposable.  
\end{corollary}

Thus, if $ \kappa, \kappa' \geq 2^ \lambda $,
then
$\lambda \stackrel{\kappa}{\Rightarrow} \mu$
if and only if 
$\lambda \stackrel{\kappa'}{\Rightarrow} \mu$.

\begin{lemma}\label{lem}
Suppose  that $ \lambda \geq \mu$  are infinite regular cardinals, 
and $\kappa  $ is a cardinal.  
Suppose that $ (f_ \beta ) _{ \beta \in \kappa } $
is a given set of functions from $ \lambda $ to $\mu$.
Then the following are equivalent.

(a) Whenever $D$ is an ultrafilter
uniform over $ \lambda $ then there exists some $ \beta \in \kappa $
such that $f_ \beta (D)$ is uniform over $ \mu$.

(b) For every function $g: \kappa \to S_\mu(\mu)$ there exists some finite
set $F \subseteq \kappa $ such that 
$ \left| \bigcap _{\beta \in F} f_\beta  ^{-1}(g( \beta )) \right| < \lambda $.
\end{lemma}

\begin{proof}
We show that the negation of (a) is equivalent to the negation of (b).

Indeed, (a) is false if and only if there exists an ultrafilter 
$D$ uniform over $\lambda $ such that for every $\beta \in \kappa $ 
it happens that $f_\beta(D)$ is not uniform over $\mu$.
This means that for every $ \beta \in \kappa $ there exists some
$ g(\beta) \in S_\mu (\mu)$ such that $g(\beta)  \in f_ \beta (D)$,
that is, $ f_\beta  ^{-1}(g( \beta )) \in D$.

Thus, there exists some $D$ which makes (a) false if and only if
there exists some function $g:\kappa\to S_\mu(\mu)$
such that the set
$ \{ f_\beta ^{-1} (g(\beta) ) | \brfr \beta \in \kappa \} \cup 
\{ x \subseteq \lambda | | \lambda \setminus x| < \lambda \}  $
 has the finite intersection property.

Equivalently, there exists some function $g:\kappa\to S_\mu(\mu)$
such that for every finite $F \subseteq \kappa $ 
the cardinality of $ \bigcap_{\beta\in F}  f_\beta ^{-1} (g(\beta)) $ 
is equal to $ \lambda  $.

This is exactly the negation of (b).
\end{proof}

Definition \ref{limpm}
suggests the introduction of the following
compactness property for logics. 

\begin{definition}\label{complog}\label{8}
Suppose that 
$ \lambda $ is an infinite cardinal, and 
$ \kappa \geq  \cf S_ \lambda  ( \lambda ) $.

Fix some set $V \subseteq S_ \lambda  ( \lambda )$
 cofinal in  $ S_ \lambda  ( \lambda )$
of cardinality $\leq \kappa $.
 
If $ (\varphi_v) _{v \in V} $ are sentences of some logic, 
and $ \gamma \in \lambda $, let
$T_ \gamma = \{ \varphi_v | v \in V, \brfr \gamma \not \in v\}$ 
and $T = \{ \varphi_v | v \in V \}$.

We say that a logic $\mathcal L$ is 
$ \kappa $-$  ( \lambda )$\emph{-compact} if and only if 
the following holds.

Whenever $ \Sigma$ is a set of sentences of $\mathcal L$,
$ |\Sigma| \leq \kappa $, $ (\varphi_v) _{v \in V} $
are sentences of $\mathcal L$ and
\[ 
 \Sigma \cup T_ \gamma   \text{ has a model for every }  \gamma \in \lambda ,
\] 
then 
\[ 
 \Sigma \cup T  \text{ has a model.}
\]
When the restriction $ |\Sigma| \leq \kappa $
on $\Sigma$ is dropped, we say that $\mathcal L$ is 
$ [ \lambda ]$\emph{-compact}. That is, 
$ [ \lambda ]$-compactness is $ \infty $-$  ( \lambda )$-compactness.
(Notice: for $ \lambda $ singular, this is a {\bf distinct} notion
from $ [ \lambda ]$-compactness as defined in \cite{Ma}.
However, the two notions coincide for $ \lambda $ regular, see
Proposition \ref{reg}).
\end{definition}

\begin{remark}\label{indeplogic}
Since $ \kappa $ is supposed to be 
$ \geq  \cf S_ \lambda  ( \lambda ) $, the above definition
is independent from the choice of the set 
$V \subseteq S_ \lambda  ( \lambda )$.

Indeed, suppose that
the logic $\mathcal L$ is 
$ \kappa $-$  ( \lambda )$-compact when the definition is
given relative to the cofinal set $V$, and let 
$V^*$ be 
 cofinal in  $ S_ \lambda  ( \lambda )$ with
$| V^*| \leq \kappa $.

Let $ (\varphi^*_w) _{w \in V^*} $ and 
$ \Sigma^*$ be sentences of $\mathcal L$, with $| \Sigma^*| \leq \kappa $,
and let $T^*_ \gamma = \{ \varphi^*_w | w \in V^*, \gamma \not \in w\}$ 
and $T^* = \{ \varphi^*_w | w \in V^* \}$.
Assuming that $ \Sigma^* \cup T^*_ \gamma  $  has a model
$ \mathfrak A_\gamma $, whenever $ \gamma \in \lambda $, we want to show that
$ \Sigma^* \cup T^* $ has a model.

For each
$ v \in V$ let us consider a new sentence $ \varphi _v$.
 The actual form of $ \varphi _v$ does not matter,
we shall simply use $ \varphi _v$ as a truth value which
 should be independent from 
all the sentences $ (\varphi^*_w) _{w \in V^*} $ and from
all sentences of $ \Sigma^*$. For example,
add a new constant symbol $c$ and new relation symbols
$ (R_v) _{v \in V} $, and let $ \varphi _v \equiv R_v(c)$.
Alternatively, let $ \varphi _v \equiv \exists x R_v(x)$.

Let $\Sigma = \Sigma^* \cup 
\{ \varphi_v \Rightarrow  \varphi^*_w | 
v \in V, w \in V^*,  v \supseteq w\}$.
Clearly, $|\Sigma|\leq \kappa $. 
As before, let
$T_ \gamma = \{ \varphi_v | v \in V, \gamma \not \in v\}$.

It is easy to show that
$\Sigma \cup T_ \gamma $ has a model for every $ \gamma \in \lambda $:
just expand the model
$ \mathfrak A_\gamma $
in such a way that $ \varphi _v$
holds true just in case $ \gamma \not \in v$.

By 
$ \kappa $-$  ( \lambda )$-compactness as given by the definition 
relative to the cofinal set $V$, we get that 
$ \Sigma \cup T$ has a model, where
$T = \{ \varphi_v | v \in V \}$.
Let $ \mathfrak A$ be a model of 
$ \Sigma \cup T$. 
$ \mathfrak A$ is also a model 
of 
$ \Sigma^* $, since
$ \Sigma \supseteq \Sigma^* $.
Moreover,
for every  $ w \in V^*$ there exists
$v \in V$ such that $ v \supseteq w$,
 since $V$ is cofinal in $ S_\lambda (\lambda )$.
Since $ \varphi _v$ holds in $ \mathfrak A$,
and 
$ \varphi_v \Rightarrow  \varphi^*_w $
is in $ \Sigma$, then $\varphi^*_w $ holds in $ \mathfrak A$. 
Thus 
$ \mathfrak A$ is a model 
of 
$ \Sigma^* \cup T^* $.

We have showed that  the logic $\mathcal L$ is 
$ \kappa $-$  ( \lambda )$-compact with reference to the definition 
given with respect to the cofinal set $ V^*$.
\end{remark}

\begin{remark}\label{9bis} 
Suppose, as in Definition \ref{complog},  that 
$ \lambda $ is an infinite cardinal, 
$ \kappa \geq  \cf S_ \lambda  ( \lambda ) $,
$V \subseteq S_ \lambda  ( \lambda )$
 is cofinal in  $ S_ \lambda  ( \lambda )$
and  $|V| \leq \kappa $.
 
If $ (\Gamma_v) _{v \in V} $ are sets of sentences of some logic $\mathcal L$, 
and $ \gamma \in \lambda $, let
$U_ \gamma =\bigcup \{ \Gamma_v | v \in V,  \gamma \not \in v\}$ 
and $U = \bigcup _{v \in V} \Gamma_v $.

Arguments similar to the ones in Remark \ref{indeplogic} 
(cf. also \cite[proof of Proposition 1.1.1]{Ma}) show that a logic $\mathcal L$ is 
$ \kappa $-$  ( \lambda )$-compact if and only if 
the following holds.

Whenever $ (\Gamma_v) _{v \in V} $ are
 sets of sentences of $\mathcal L$, and
$|\Gamma_v| \leq \kappa $ for every $v \in V$,
if $U_ \gamma $ has a model for every $ \gamma \in \lambda $,
then $U$ has a model.
\end{remark}

\begin{proposition}\label{cfcp}
If $\kappa \geq  \cf S_ \lambda  ( \lambda ) $, then every
$ \kappa $-$  ( \lambda )$-compact logic
is $ \kappa $-$  ( \cf \lambda )$-compact.
 \end{proposition}

Recall the definition 
of $ \kappa $-$ ( \lambda, \lambda  )$-compactness from
\cite{easter} or 
\cite[Proposition 2.1]{bumi}. Cf. also \cite[Section 1]{arxiv}.
Recall the definition 
of $ [ \mu, \lambda  ]$-compactness from
\cite[Definition 1.1.3]{Ma}. 

\begin{proposition}\label{reg}
If $ \lambda $ is an infinite regular cardinal, and 
$ \kappa \geq  \lambda   $, then a logic  is 
$ \kappa $-$ ( \lambda )$-compact if and only if it is
$ \kappa $-$ ( \lambda, \lambda  )$-compact. 

If $ \lambda $ is an infinite regular cardinal, then
a logic is $[ \lambda ]$-compact if and only if it is
$ [ \lambda, \lambda  ]$-compact. 
\end{proposition}

Notice that in Remarks \ref{indeplogic} and \ref{9bis},  and in Propositions 
\ref{cfcp} and \ref{reg}
we need really few regularity properties for $\mathcal L$.
In contrast, the proofs of the next propositions make heavy use of relativization.
See \cite{Eb} for the definition of a regular logic.
In the next propositions slightly less is needed:
instead of requiring a logic $\mathcal L$ to be regular, it is enough to require that
$\mathcal L$ satisfies the properties listed in \cite[Section 1]{jsl}.

\begin{proposition}\label{impimpcomp}
Suppose that 
$ \lambda \geq \mu$ are infinite cardinals, and 
$ \kappa \geq \sup \{ \cf S_ \lambda  ( \lambda ), \cf S_\mu(\mu)\} $.

If $\lambda \stackrel{\kappa}{\Rightarrow} \mu$
holds, then every
$ \kappa $-$ ( \lambda )$-compact regular logic is
$ \kappa $-$ ( \mu )$-compact.

If every ultrafilter uniform over $ \lambda $ is 
$\mu$-decomposable, then every
$ [ \lambda ]$-compact regular logic is
$ [ \mu ]$-compact.
\end{proposition}

\begin{proof}
Similar to the proof of \cite[Theorem 2.4]{bumi}. Cf. also \cite[Proposition 1]{easter}. For the last statement, use Corollary \ref{coruf}. 
\end{proof}
 
Recall from  \cite[p.440]{jsl} the definition a limit ultrafilter.
We say that a limit ultrafilter $(D,V)$ (where $D$ is over $I$) is $ \lambda $-\emph{decomposable}
if and only if there exists a function $f:I\to \lambda $ such that 
$ \{ (i,j)|f(i)=f(j)\}\in V $ and for every $v \subseteq \lambda $ with
$|v|<\lambda $ we have $f^{-1}(v)\not \in D$.

Recall from  \cite[Definition 2.4]{jsl} the definition of the family $ \mathbf{F}(\mathcal L) $ associated to a logic 
$\mathcal L $.

\begin{proposition}\label{abstractjsl}
If $\mathcal L$ is a regular logic, then the following are equivalent.

(a) $\mathcal L$ is $ [ \lambda ]$\emph-compact.

(b) For every relativized (that is, many sorted) expansion 
${\mathfrak A}$ of the model 
$ \mathfrak A'$
introduced in Definition \ref{limpm}, 
there exist an $\mathcal L$-elementary extension 
$\mathfrak B $
of
$\mathfrak A$ and an element $x \in B$ such that 
$ {\mathfrak B} \models \neg R_v( x) $ for every $ v \in V$
(the choice of the set $V$ cofinal in $S_\lambda (\lambda )$ is
not relevant.)

(c) For every infinite cardinal $\nu$ there exists
 some $ \lambda $-decomposable  limit ultrafilter
$(D,V)$ such that $(D,V,\nu) \in \mathbf{F}(\mathcal L)$.
 \end{proposition}

\begin{proof}
Similar to the proof of \cite[Theorem 2.3]{jsl}, using Remark \ref{6x}. 
\end{proof}

\end{document}